\documentclass{gtart}

\def\ifplaintex{\expandafter\ifx\csname documentclass\endcsname\relax}

\def\gtp{{\mathsurround=0pt\it $\cal G\mskip-2mu$eometry \&\ 
$\cal T\!\!$opology $\cal P\!$ublications}}  

\def\recd{{\small Received:\qua\receiveddate\ifx\reviseddate\relax
\else\qquad Revised:\qua\reviseddate\fi\par}} 

\def\lognumber#1{\def\thelognumber{#1}}
\def\volumenumber#1{\def\thevolumenumber{#1}}
\def\volumeyear#1{\def\thevolumeyear{#1}}
\def\papernumber#1{\def\thepapernumber{#1}}
\def\pagenumbers#1#2{\def\startpage{#1}\def\finishpage{#2}}
\def\published#1{\def\publishdate{#1}}

\def\received#1{\def\receiveddate{#1}}
\def\revised#1{\def\reviseddate{#1}}
\def\accepted#1{\def\accepteddate{#1}}

\def\asciiaddress#1{\def\theasciiaddress{#1}}

\long\def\asciiabstract#1{\long\def\theasciiabstract{#1}}

\let\\\par\let\thelognumber\relax\let\thevolumenumber\relax
\let\thepapernumber\relax\let\thevolumeyear\relax\let\startpage\relax
\let\finishpage\relax\let\publishdate\relax\let\receiveddate\relax
\let\reviseddate\relax\let\accepteddate\relax\let\theasciititle\relax
\let\theasciiauthors\relax\let\theasciiaddress\relax
\let\theasciiabstract\relax

\let\theasciiemail\relax

\ifplaintex
\font\logobig=cmssbx10 scaled 3836
\font\logomed=cmssbx10 scaled 2557
\else
\font\logobig=cmssbx10 scaled 4200
\font\logomed=cmssbx10 scaled 2800
\fi

\long\def\makeagttitle{   
\count0=\startpage
\agt\hfill      
\hbox to 45truept{\vbox to 0pt{\vglue -13truept{\logomed A\kern -.37em{\logobig 
T}\kern -.38em G}\vss}\hss}
\break
{\small Volume \thevolumenumber\ (\thevolumeyear)
\startpage--\finishpage\nl
Published: \publishdate}

\vglue .25truein

{\parskip=0pt\leftskip 0pt plus
1fil\def\\{\par\smallskip}{\Large\bf\thetitle}\par\medskip} \vglue
0.05truein

{\parskip=0pt\leftskip 0pt plus 1fil\def\\{\par}{\sc\theauthors}
\par\medskip}%
 
\vglue 0.03truein 

{\small\leftskip 25truept\rightskip 25truept{\bf Abstract}\stdspace\theabstract

{\bf AMS Classification}\stdspace\theprimaryclass
\ifx\thesecondaryclass\relax\else; \thesecondaryclass\fi\par
{\bf Keywords}\stdspace \thekeywords\par}\vglue 7truept

}   

\ifplaintex
\font\phead=cmsl9 scaled 950
\font\pnum=cmbx10 scaled 913
\font\pfoot=cmsl9 scaled 950
\headline{\vbox to 0pt{\vskip -4.5mm\line{\small\phead\ifnum
\count0=\startpage ISSN 1472-2739 (on-line) 1472-2747 (printed)
\hfill {\pnum\folio}\else\ifodd\count0\def\\{ }%
\ifx\theshorttitle\relax\thetitle\else\theshorttitle\fi\hfill{\pnum\folio}
\else\def\\{ and }{\pnum\folio}\hfill\ifx\theshortauthors\relax\theauthors
\else\theshortauthors\fi\fi\fi}\vss}}
\footline{\vbox to 0pt{\vglue 0mm\line{\small\pfoot\ifnum\count0=\startpage
\copyright\ \gtp\hfill\else
\agt, Volume \thevolumenumber\ (\thevolumeyear)\hfill\fi}\vss}}
\else
\font=cmsl9 scaled 1050
\font\lnum=cmbx10 
\font=cmsl9 scaled 1050
\makeatletter
\def\@oddhead{{\small\lhead\ifnum\count0=\startpage ISSN 1472-2739 
(on-line) 1472-2747 (printed)\hfill {\lnum\number\count0}\else\ifodd\count0
\def\\{ }\ifx\theshorttitle\relax \thetitle \else\theshorttitle\fi\hfill
{\lnum\number\count0}\else\def\\{ and }{\lnum\number\count0}
\hfill\ifx\theshortauthors\relax 
\theauthors\else\theshortauthors\fi\fi\fi}}\def\@evenhead{\@oddhead}
\def\@oddfoot{\small\lfoot\ifnum\count0=\startpage\copyright\ \gtp\hfill\else
\agt, Volume \thevolumenumber\ (\thevolumeyear)\hfill\fi}
\def\@evenfoot{\@oddfoot}
\makeatother
\fi
\let\maketitlepage\makeagttitle
\let\makeshorttitle\maketitlepage
\let\maketitle\maketitlepage

\newwrite\gtoutfile
\long\gdef\makeheadfile{  
{\def\\{, }\def\s{ }
\immediate\openout\gtoutfile head.xxx
\immediate\write\gtoutfile{To: math@arxiv.org}
\immediate\write\gtoutfile{Subject: put OR rep NNNNN:ppppp}
\immediate\write\gtoutfile{--text follows this line--}
\immediate\write\gtoutfile{Proxy-for: \ifx\theasciiauthors\relax
\theauthors\else\theasciiauthors\fi\s<\ifx\theasciiemail\relax\theemail\else\theasciiemail\fi>}
\immediate\write\gtoutfile{\noexpand\\}
\immediate\write\gtoutfile{Authors: \ifx\theasciiauthors\relax
\theauthors\else\theasciiauthors\fi}
{\def\\{ }\immediate\write\gtoutfile{Title: \ifx\theasciititle\relax
\thetitle\else\theasciititle\fi}}
\immediate\write\gtoutfile{Subj-class: GT or SG, GR etc}
\immediate\write\gtoutfile{MSC-class: \theprimaryclass\ifx\thesecondaryclass\relax\else, \thesecondaryclass\fi}
\immediate\write\gtoutfile{Journal-ref: Algebr. Geom. Topol. \thevolumenumber\s
(\thevolumeyear) \startpage-\finishpage}
\immediate\write\gtoutfile{Comments: Published by Algebraic and
Geometric Topology at}
\immediate\write\gtoutfile{\s\s\s  http://www.maths.warwick.ac.uk/agt/AGTVol\thevolumenumber/agt-\thevolumenumber-\thepapernumber.abs.html}
\immediate\write\gtoutfile{\noexpand\\}
\immediate\write\gtoutfile{}
\ifx\theasciiabstract\relax
\immediate\write\gtoutfile{\theabstract}\else
\immediate\write\gtoutfile{\theasciiabstract}\fi
\immediate\write\gtoutfile{}
\immediate\write\gtoutfile{\noexpand\\}
\immediate\write\gtoutfile{}
\immediate\closeout\gtoutfile}}  

\def\maketitlepage{\makeagttitle\makeheadfile}
\let\makeshorttitle\maketitlepage
\let\maketitle\maketitlepage

\def\ifplaintex{\expandafter\ifx\csname documentclass\endcsname\relax}

\def\gtp{{\mathsurround=0pt\it $\cal G\mskip-2mu$eometry \&\ 
$\cal T\!\!$opology $\cal P\!$ublications}}  

\def\recd{{\small Received:\qua\receiveddate\ifx\reviseddate\relax
\else\qquad Revised:\qua\reviseddate\fi\par}} 

\def\lognumber#1{\def\thelognumber{#1}}
\def\volumenumber#1{\def\thevolumenumber{#1}}
\def\volumeyear#1{\def\thevolumeyear{#1}}
\def\papernumber#1{\def\thepapernumber{#1}}
\def\pagenumbers#1#2{\def\startpage{#1}\def\finishpage{#2}}
\def\published#1{\def\publishdate{#1}}

\def\received#1{\def\receiveddate{#1}}
\def\revised#1{\def\reviseddate{#1}}
\def\accepted#1{\def\accepteddate{#1}}

\def\asciiaddress#1{\def\theasciiaddress{#1}}

\long\def\asciiabstract#1{\long\def\theasciiabstract{#1}}

\let\\\par\let\thelognumber\relax\let\thevolumenumber\relax
\let\thepapernumber\relax\let\thevolumeyear\relax\let\startpage\relax
\let\finishpage\relax\let\publishdate\relax\let\receiveddate\relax
\let\reviseddate\relax\let\accepteddate\relax\let\theasciititle\relax
\let\theasciiauthors\relax\let\theasciiaddress\relax
\let\theasciiabstract\relax

\let\theasciiemail\relax

\ifplaintex
\font\logobig=cmssbx10 scaled 3836
\font\logomed=cmssbx10 scaled 2557
\else
\font\logobig=cmssbx10 scaled 4200
\font\logomed=cmssbx10 scaled 2800
\fi

\long\def\makeagttitle{   
\count0=\startpage
\agt\hfill      
\hbox to 45truept{\vbox to 0pt{\vglue -13truept{\logomed A\kern -.37em{\logobig 
T}\kern -.38em G}\vss}\hss}
\break
{\small Volume \thevolumenumber\ (\thevolumeyear)
\startpage--\finishpage\nl
Published: \publishdate}

\vglue .25truein

{\parskip=0pt\leftskip 0pt plus
1fil\def\\{\par\smallskip}{\Large\bf\thetitle}\par\medskip} \vglue
0.05truein

{\parskip=0pt\leftskip 0pt plus 1fil\def\\{\par}{\sc\theauthors}
\par\medskip}%
 
\vglue 0.03truein 

{\small\leftskip 25truept\rightskip 25truept{\bf Abstract}\stdspace\theabstract

{\bf AMS Classification}\stdspace\theprimaryclass
\ifx\thesecondaryclass\relax\else; \thesecondaryclass\fi\par
{\bf Keywords}\stdspace \thekeywords\par}\vglue 7truept

}   

\ifplaintex
\font\phead=cmsl9 scaled 950
\font\pnum=cmbx10 scaled 913
\font\pfoot=cmsl9 scaled 950
\headline{\vbox to 0pt{\vskip -4.5mm\line{\small\phead\ifnum
\count0=\startpage ISSN 1472-2739 (on-line) 1472-2747 (printed)
\hfill {\pnum\folio}\else\ifodd\count0\def\\{ }%
\ifx\theshorttitle\relax\thetitle\else\theshorttitle\fi\hfill{\pnum\folio}
\else\def\\{ and }{\pnum\folio}\hfill\ifx\theshortauthors\relax\theauthors
\else\theshortauthors\fi\fi\fi}\vss}}
\footline{\vbox to 0pt{\vglue 0mm\line{\small\pfoot\ifnum\count0=\startpage
\copyright\ \gtp\hfill\else
\agt, Volume \thevolumenumber\ (\thevolumeyear)\hfill\fi}\vss}}
\else
\font=cmsl9 scaled 1050
\font\lnum=cmbx10 
\font=cmsl9 scaled 1050
\makeatletter
\def\@oddhead{{\small\lhead\ifnum\count0=\startpage ISSN 1472-2739 
(on-line) 1472-2747 (printed)\hfill {\lnum\number\count0}\else\ifodd\count0
\def\\{ }\ifx\theshorttitle\relax \thetitle \else\theshorttitle\fi\hfill
{\lnum\number\count0}\else\def\\{ and }{\lnum\number\count0}
\hfill\ifx\theshortauthors\relax 
\theauthors\else\theshortauthors\fi\fi\fi}}\def\@evenhead{\@oddhead}
\def\@oddfoot{\small\lfoot\ifnum\count0=\startpage\copyright\ \gtp\hfill\else
\agt, Volume \thevolumenumber\ (\thevolumeyear)\hfill\fi}
\def\@evenfoot{\@oddfoot}
\makeatother
\fi
\let\maketitlepage\makeagttitle
\let\makeshorttitle\maketitlepage
\let\maketitle\maketitlepage

\newwrite\gtoutfile
\long\gdef\makeheadfile{  
{\def\\{, }\def\s{ }
\immediate\openout\gtoutfile head.xxx
\immediate\write\gtoutfile{To: math@arxiv.org}
\immediate\write\gtoutfile{Subject: put OR rep NNNNN:ppppp}
\immediate\write\gtoutfile{--text follows this line--}
\immediate\write\gtoutfile{Proxy-for: \ifx\theasciiauthors\relax
\theauthors\else\theasciiauthors\fi\s<\ifx\theasciiemail\relax\theemail\else\theasciiemail\fi>}
\immediate\write\gtoutfile{\noexpand\\}
\immediate\write\gtoutfile{Authors: \ifx\theasciiauthors\relax
\theauthors\else\theasciiauthors\fi}
{\def\\{ }\immediate\write\gtoutfile{Title: \ifx\theasciititle\relax
\thetitle\else\theasciititle\fi}}
\immediate\write\gtoutfile{Subj-class: GT or SG, GR etc}
\immediate\write\gtoutfile{MSC-class: \theprimaryclass\ifx\thesecondaryclass\relax\else, \thesecondaryclass\fi}
\immediate\write\gtoutfile{Journal-ref: Algebr. Geom. Topol. \thevolumenumber\s
(\thevolumeyear) \startpage-\finishpage}
\immediate\write\gtoutfile{Comments: Published by Algebraic and
Geometric Topology at}
\immediate\write\gtoutfile{\s\s\s  http://www.maths.warwick.ac.uk/agt/AGTVol\thevolumenumber/agt-\thevolumenumber-\thepapernumber.abs.html}
\immediate\write\gtoutfile{\noexpand\\}
\immediate\write\gtoutfile{}
\ifx\theasciiabstract\relax
\immediate\write\gtoutfile{\theabstract}\else
\immediate\write\gtoutfile{\theasciiabstract}\fi
\immediate\write\gtoutfile{}
\immediate\write\gtoutfile{\noexpand\\}
\immediate\write\gtoutfile{}
\immediate\closeout\gtoutfile}}  

\def\maketitlepage{\makeagttitle\makeheadfile}
\let\makeshorttitle\maketitlepage
\let\maketitle\maketitlepage

\lognumber{23}
\volumenumber{1}
\volumeyear{2001}
\papernumber{23}
\pagenumbers{445}{468}
\received{24 March 2001}
\revised{17 August 2001}
\accepted{17 August 2001}
\published{23 August 2001}

\usepackage{amssymb,graphicx}
\let\endpf\endproof

 \newtheorem{defn}{Definition}[section]
 \newtheorem{lemma}[defn]{Lemma}
 
 \newtheorem{theorem}[defn]{Theorem}
 \newtheorem{definition}[defn]{Definition}
 \newtheorem{remark}[defn]{Remark}
 \newtheorem{prop}[defn]{Proposition}
 \newtheorem{corollary}[defn]{Corollary}

 \newtheorem{claim}[defn]{Claim}

 \newcommand{\zz}{{\bf Z}}

 \newcommand{\calc}{{\cal C}}

 \newcommand{\cala}{{\cal A}}

 \newcommand{\calh}{{\cal H}}

 \newcommand{\calg}{{\cal G}}

 \newcommand{\calm}{{\cal M}}

 \newcommand{\autftwo}{\textrm{Aut}(F_2)}
 \newcommand{\outftwo}{\textrm{Out}(F_2)}
 \newcommand{\autfn}{\textrm{Aut}(F_n)}
 \newcommand{\outfn}{\textrm{Out}(F_n)}

 \newcommand{\autg}{\textrm{Aut}(G)}
 
 \newcommand{\mgzz}{\calm_{g,0,0}}
 \newcommand{\mgoz}{\calm_{g,1,0}}
 \newcommand{\mgzo}{\calm_{g,0,1}}
 \newcommand{\mtzz}{\calm_{2,0,0}}
 \newcommand{\mgbn}{\calm_{g,b,n}}
 
 \newcommand{\mgplusk}{\calm_{g+k,0,0}}
 
 \newcommand{\sgbn}{S_{g,b,n}}
 \newcommand{\sgzz}{S_{g,0,0}}
 
 \newcommand{\sgzo}{S_{g,0,1}}

 \newcommand{\pigzz}{\pi_1 (\sgzz)}
 \newcommand{\pisg}{\pi_1 (\sgzz)}
 \newcommand{\pigzo}{\pi_1 (\sgzo)}
 \newcommand{\autpisg}{\textrm{Aut} ( \pi_1 (S_{g,0,0}))}
 \newcommand{\outpisg}{\textrm{Out} ( \pi_1 (S_{g,0,0}))}
 \newcommand{\innpisg}{\textrm{Inn} ( \pi_1 (S_{g,0,0}))}

 \newcommand{\inv}{{^{-1}}}
 \newcommand{\calb}{{\cal B}}
 \newcommand{\g}{{\gamma}}
 \newcommand{\G}{{\Gamma}}
 \newcommand{\e}{{\epsilon}}
 \newcommand{\st}{{{\textrm{\it{start}}}}}
 \newcommand{\en}{{{\textrm{\it{end}}}}}
 \newcommand{\al}{{\alpha}}
 \newcommand{\lan}{{\langle}}
 \newcommand{\ran}{{\rangle}}
 \newcommand{\s}{{\sigma}}

\begin{document}

 \title{On the linearity problem for mapping class groups}
 \authors{Tara E. Brendle\\Hessam Hamidi-Tehrani}
\email{tbrendle@math.columbia.edu, hessam@math.columbia.edu}
\address{Columbia University\\Department of Mathematics\\
New York, NY  10027, USA}
\secondaddress{B.C.C. of the City University of New York\\Department of Mathematics and Computer Science\\Bronx, NY  10453, USA}
\asciiaddress{Columbia University\\Department of Mathematics\\
New York, NY  10027, USA\\ 
B.C.C. of the City University of New York\\Department of Mathematics and Computer Science\\Bronx, NY  10453, USA}

 \begin{abstract}

Formanek and Procesi have demonstrated that $\autfn$ is not linear for
 $n \geq 3$.  Their technique is to construct nonlinear groups of a
 special form, which we call {\it FP-groups}, and then to embed a
 special type of automorphism group, which we call a {\it poison
 group}, in $\autfn$, from which they build an FP-group.  We first
 prove that poison groups cannot be embedded in certain mapping class
 groups. We then show that no FP-groups of any form can be embedded in
 mapping class groups.  Thus the methods of Formanek and Procesi fail
 in the case of mapping class groups, providing strong evidence that
 mapping class groups may in fact be linear.  

\end{abstract}
\asciiabstract{
Formanek and Procesi have demonstrated that Aut(F_n) is not linear for
n >2.  Their technique is to construct nonlinear groups of a
special form, which we call FP-groups, and then to embed a
special type of automorphism group, which we call a poison
group, in Aut(F_n), from which they build an FP-group.  We first
prove that poison groups cannot be embedded in certain mapping class
groups. We then show that no FP-groups of any form can be embedded in
mapping class groups.  Thus the methods of Formanek and Procesi fail
in the case of mapping class groups, providing strong evidence that
mapping class groups may in fact be linear.  
}
\maketitle

 \section{Introduction}
 The question of whether mapping class groups are linear has been around
 for some time.  The recent work of Bigelow \cite{bi} and also Krammer
 \cite{kr} in  determining that the braid group is linear has renewed
 interest in the subject, due to the close relationship between mapping
 class groups and braid groups. Let $\sgbn$ denote a surface of genus
 $g$ with $b$ boundary components and $n$ fixed points.  Let $\mgbn$ denote
 the mapping class group of $\sgbn$.  We assume throughout that maps fix boundary components pointwise.  Bigelow and Budney \cite{bb} and
 independently Korkmaz \cite{ko} recently determined that $\mtzz$ is linear.
  Korkmaz also showed in \cite{ko} that mapping class groups contain very
 large linear subgroups, namely, the hyperelliptic subgroups.  However,
 the question of linearity remains open for mapping class groups of surfaces
 of genus 3 or greater.   

 Let $F_n$ denote the free group of rank $n$.  It is well known that $\outftwo$ and $\autftwo$ are linear.  The former fact is due to Nielsen \cite{n}, and the latter follows by \cite{dfg} from the linearity of the 4-string braid group $B_4$, which is due to Krammer \cite{kr}.

On the other hand, Formanek and Procesi demonstrated in \cite{fp} that
 $\autfn$ is not a linear group for $n \geq 3$.  A simple corollary of this result is that $\outfn$
 is not linear for $n \geq 4$.  The well-known fact due to Nielsen \cite{mks}
 that $\mgzz$ is isomorphic to $\outpisg$ suggests that it may be possible
 to apply the methods of Formanek and Procesi to mapping class groups,
 though it may not be immediately clear how to do so.  

 Formanek and Procesi define a class of nonlinear groups, which we will generalize slightly and refer to as {\it Formanek and Procesi groups}, or {\it FP-groups} for short.   We will show that the existence of FP-subgroups of $\mgzo$ would imply that $\mgplusk$ is not linear for $k \geq 1$.  We will also focus our attention on a special kind of automorphism group, which we call a {\it poison group}.  We will describe the particular method of Formanek and Procesi for constructing FP-groups from poison subgroups.

 This work originated in an attempt to use the methods of Formanek and Procesi to show that $\mgzz$ is not linear for $g \geq 3$.  We prove instead that the essential building blocks of the Formanek and Procesi method do not exist in mapping class groups, first in a special case.  

\medskip{\bf Theorem A}\qua{\sl Poison subgroups cannot be embedded in $\mgzo$. }

\medskip

Thus the particular technique of Formanek and Procesi fails to show that certain mapping class groups are not linear.  We then generalize this result as follows.

\medskip

{\bf Theorem B}\qua{\sl FP-groups do not embed in $\mgbn$ for any $g$, $b$, and $n$.}

\medskip

  Our paper is organized as follows.  
 In Section 2, we give an overview of the methods of Formanek and Procesi
 for constructing a nonlinear subgroup of $\autfn$ from a poison subgroup.
  In Section 3, we establish connections between certain mapping class
 groups and the automorphism group of a closed surface.  In Section 4
 we prove Theorem A.  In Section 5 we prove Theorem B using very different techniques from those used in Section 4.  Though Theorem A is a special case of Theorem B, we include a separate proof of Theorem A both for the sake of highlighting the particular construction of Formanek and Procesi and also because the methods used are interesting in their own right.  The reader should note, however, that Sections 3, 4, and 5 are completely independent of one another.  For example, the reader interested only in Theorem B could read Sections 1, 2, and 5 without any loss of continuity.

\noindent {\bf Acknowledgements}\qua
 The authors would like to express their sincere gratitude to Joan Birman
 and Alex Lubotzky for suggesting the search for poison subgroups in mapping
 class groups, and also to Matthew Zinno for helping to point us in the
 other direction.  We
 thank all three, as well as Walter Neumann, Brian Mangum, Gabriel Rosenberg,
 and Abhijit Champanerkar for many useful discussions.  We are also grateful to the referee for many helpful questions and suggestions.

The first author was partially supported under NSF Grant DMS-9973232. The second author was partially supported by PSC-CUNY Research Grant 63463 00 32.

 \section{The method of Formanek and Procesi}
 Let $G$ be any group, and let $\calh (G)$ denote the following HNN-extension
 of $G \times G$:
 $$
 \calh (G) = \langle G \times G, t \hskip .1in | \hskip .1in t (g, g)
 t^{-1} = (1, g), g \in G  \rangle .
 $$
 In other words, conjugation by $t$ in the HNN-extension carries the diagonal
 subgroup $G \times G$ onto its second factor.  Formanek and Procesi show
 in the following theorem that such groups exhibit special behavior under
 a linear representation.
   
 \begin{theorem}[Formanek and Procesi, \cite{fp}]\label{thmfour}
 Let $G$ be a group.
 Then the image of the subgroup $G \times \{ 1 \}$ under any linear representation
 of $\calh (G)$ is nilpot\-ent-by-abelian-by-finite. 
 \end{theorem}

 \begin{corollary}\label{cor}
 Let $G$ be a group, and $K$ a normal subgroup of $\calh (G)$ such that the image of $G \times \{1\}$ in $\calh (G) / K$ is not nilpotent-by-abelian-by-finite.  Then $\calh (G) / K$ is not linear.
 \end{corollary}

\proof
Let $\rho : \calh(G)/K \rightarrow GL_N (k)$ be a linear representation where $k$ is a field.  Let $\pi : \calh (G) \rightarrow \calh(G)/K$ be the natural projection map.  Then $\rho \circ \pi$ is a linear representation of $\calh(G)$ and hence by Theorem~\ref{thmfour}, $\rho(\pi(G \times \{1\}))$ is nilpotent-by-abelian-by-finite.  Thus $\rho$ is not faithful.  
\endpf

 We will call a group of the type described in Corollary~\ref{cor}
 a
 {\it Formanek and Procesi group}, or {\it FP-group} for short.  
  We now describe the particular construction of Formanek and Procesi in demonstrating the nonlinearity of $\autfn$ for $n \geq 3$.

 Let $G$ be any group.  Let $x_1, x_2, x_3$ be elements of $G$ such that
 $\langle x_1, x_2, x_3 \rangle \cong F_3$.  Let $\phi_1, \phi_2 \in \autg$
 be two maps such that
 \begin{enumerate}
 \item $\phi_i (x_j) = x_j,$ \hskip .15in $ i, j = 1, 2$, \hskip .15in
 and
 \item $\phi_i (x_3) = x_3 x_i,$ \hskip .15in $i = 1, 2$.
 \end{enumerate}
 We will call the subgroup $\langle \phi_1, \phi_2 \rangle$ a {\it poison
 subgroup} of $\autg$.  We can define poison subgroups of the mapping class group $\mgzo$ analogously, since in this case the mapping class group acts on $\pigzo$.  Notice that the second condition implies that $\langle \phi_1,
 \phi_2 \rangle \cong F_2$.  Thus poison groups, being isomorphic to the linear group $F_2$, are not themselves a kind of FP-group.  However, as the following lemma shows, their existence in an automorphism group $\autg$ implies that $\autg$ is not linear (hence the name ``poison groups'', though it suggests a bias towards linearity).

 \begin{lemma}\label{hessamlemma}
Let $G$ be any group.  If $\autg$ contains a poison subgroup, then it contains an FP-subgroup isomorphic to $\calh(F_2)$.
 \end{lemma}

 \proof
 Let $\langle \phi_1, \phi_2 \rangle$ be a poison subgroup in $\autg$.  
 Following Formanek and Procesi's argument in \cite{fp}, let $\alpha_i \in \autg$ denote conjugation by $x_i$.  Consider the
 group 
 $$H = \hskip .1in \langle \phi_1, \phi_2, \alpha_1, \alpha_2, \alpha_3
 \rangle .$$  First, note that $\langle  \alpha_1, \alpha_2, \alpha_3
 \rangle$ is a normal subgroup of $H$ since both $\phi_1$ and $\phi_2$
 preserve the subgroup $\langle x_1, x_2, x_3 \rangle $.  Now let $w(a,b)$
 denote any non-trivial reduced word in the free group on the letters $a$ and $b$. By definition
 of a poison subgroup, we know that $w(\phi_1, \phi_2)(x_i) = x_i$ for i = 1,2.  This tells us that if $w(\phi_1,\phi_2)$ is in $\langle \alpha_1, \alpha_2, \alpha_3 \rangle$, then $w(\phi_1, \phi_2)$ must induce conjugation by an element in $\langle x_1,x_2,x_3 \rangle \cong F_3$, which commutes with $x_1$ and $x_2$.  But the only such element is the identity.  Hence $w (\phi_1, \phi_2)$ must be the identity map.  But we know this is not the case since
 \begin{equation}\label{word}
 w(\phi_1, \phi_2) (x_3) = x_3 w(x_1, x_2).  
 \end{equation}
 This tells us that the images of $\phi_1$ and $\phi_2$ mod $\langle 
 \alpha_1, \alpha_2, \alpha_3 \rangle $ will generate a free group.  Clearly,
 the images of $\phi_1$ and $\phi_2$ also generate the quotient of $H$
 by $\langle  \alpha_1, \alpha_2, \alpha_3 \rangle$ , and so we have a
 split exact sequence
 \begin{equation}
 1 \rightarrow \langle  \alpha_1, \alpha_2, \alpha_3 \rangle \rightarrow
 H \rightarrow \langle \phi_1 , \phi_2 \rangle \rightarrow 1. 
 \end{equation}
 Thus the only relations we have in a presentation for $H$ are given by
 conjugation, as follows:
  \begin{equation}\label{pres}
 H = \langle \phi_1, \phi_2, \alpha_1, \alpha_2, \alpha_3 \hskip .1in
 | \hskip .1in \phi_i \alpha_j \phi_i^{-1} = \alpha_j, \hskip .1in \phi_i
 \alpha_3 \phi_i^{-1} = \alpha_3 \alpha_i, \hskip .1in i, j = 1, 2 \rangle
 .
 \end{equation} 
 Rewriting the second set of relations, we obtain $\alpha_3 (\alpha_i
 \phi_i) \alpha_3^{-1} = \phi_i, \hskip .15in i=1,2$.  Since $\langle
  \phi_1, \phi_2 \rangle  \cong \langle  \alpha_1, \alpha_2 \rangle \cong
 F_2$, we have that $H \cong \calh(F_2)$, with $\alpha_3$ playing the
 role of the element $t$.  Since $F_2$ is not nilpotent-by-abelian-by-finite,
 $\calh (F_2)$ is an FP-group.  \endpf

 \section{The connection with mapping class groups}
 Our motivation for the work in this paper is the following observation,
 the proof of which we defer to the end of the section.

 \begin{claim} \label{keyclaim}
 If a poison subgroup exists in $\mgzo$ for $g \geq 2$, then the groups
 $\mgplusk$ are not linear for $k \geq 1$.    
 \end{claim}

 We have been abusing terminology a bit by talking about poison subgroups
 in $\mgzo$ and also in the context of automorphism groups.
 The distinction between the two contexts is unnecessary for our purposes,
 as the following lemma shows, since these mapping class groups are isomorphic
 to automorphism groups.

 \begin{lemma}\label{isolemma}
 $\mgzo \cong \autpisg,$ for $g \geq 2.$
 \end{lemma}

 \proof
 We begin with the exact sequence
 $$
 1 \rightarrow \innpisg \rightarrow \autpisg \rightarrow \outpisg \rightarrow
 1.
 $$
By the well-known theorem of Nielsen \cite{mks}, we have that $\outpisg \cong
 \mgzz$.  In addition, since $\pisg$ is centerless,  we can replace $\innpisg$
 with $\pigzz$ (see, for example, \cite{br}) to obtain
 \begin{equation}
 1 \rightarrow \pisg \rightarrow \autpisg \rightarrow \mgzz \rightarrow
 1.
 \label{short1}
 \end{equation}
 By \cite{bir}, we also have the following exact sequence:
 \begin{equation}
 1 \rightarrow \pisg \rightarrow \mgzo \rightarrow \mgzz \rightarrow 1.
 \label{short2}
 \end{equation}

Every short exact sequence $1 \rightarrow N \rightarrow E \rightarrow G \rightarrow 1$ induces a homomorphism $G \rightarrow \textrm{Out}(N)$, defined as follows.  Let $g \in G$, and let $e_g$ be a lift of $g \in E$.  Now, $E$ acts on $N$ by conjugation, hence we can think of $e_g$ as an element of $\textrm{Aut}(N)$.  However, since $N$ is not necessarily abelian, this map is only well defined up to conjugation by an element of $N$.  Thus we get a map $G \rightarrow \textrm{Out}(N)$.  
According to Corollary 6.8 of \cite{br}, given any short exact sequence as above, with $N$ centerless, there is a unique ``middle group" $E$ corresponding to any given homomorphism $G \rightarrow \textrm{Out}(N)$.

In Sequence~\ref{short1} above, it is clear that the map induced is the Nielsen isomorphism between $\mgzz$ and $\textrm{Out}(\pisg)$.  
In Sequence~\ref{short2}, as discussed in \cite{bir}, the image of a generator $a$ of $\pisg$ is the so-called ``spin map" associated to each curve, which induces conjugation by that curve, but can be more easily understood as a product of opposite Dehn twists about the boundary of an annular neighborhood of the curve $a$.  In other words, if $\alpha$ and $\beta$ are the two boundary curves, then the spin map associated to the curve $a$ can be written as $T_{\alpha} T_{\beta}^{-1}$, where $T_{\gamma}$ denotes the Dehn twist about the curve $\gamma$.  Let $\phi \in \mgzz$, and let $\tilde \phi$ denote a lift of $\phi$ in $\mgzo$.  Then $\tilde \phi T_{\alpha} T_{\beta}^{-1} \tilde \phi^{-1} = T_{\tilde \phi(\alpha)}T_{\tilde \phi(\beta)}^{-1}$, which is precisely the spin map associated to $\tilde \phi(a)$.  Thus, we are simply looking at the action of $\tilde \phi$ on $\pisg$, but since $\phi$ does not necessarily fix the basepoint, $\phi$ is getting mapped to the class of $\tilde \phi$ in $\textrm{Aut}$, modulo inner automorphisms.  In other words, the induced map from $\mgzz \rightarrow \textrm{Out}(\pisg)$ is also the Nielsen isomorphism.  
Now since $\pisg$ has a trivial center, we apply Corollary 6.8 of \cite{br}, and the lemma is proved.  \endpf

 \begin{remark}
{\rm The isomorphism given in Lemma~\ref{isolemma} has received some attention
 in the literature, though perhaps not as much as it deserves.
 The map itself is the obvious one, namely, any homeomorphism of a surface
 with one fixed point induces a natural automorphism of the fundamental
 group of the closed surface with the fixed point taken as base point.
  From the geometric point of view, it is not immediately clear that this
 map from $\mgzo$ to $\autpisg$ should be a surjection, i.e., it is not
 necessarily obvious that all elements of $\autpisg$ should be topologically
 induced.}  
 \end{remark}

 \begin{lemma}\label{alexlemma}
 If $\autpisg$ is not linear, then $\mgoz$ is not linear.
 \end{lemma}

 Before proving the lemma, we make
 a few observations.  From Chapter 4, Section 1 of \cite{bir} and Lemma~\ref{isolemma} we have the short exact sequence
 \begin{equation}\label{surjection}
 1 \rightarrow \zz \rightarrow \mgoz \rightarrow \autpisg \rightarrow
 1.
 \end{equation}
We note that $\zz$ is actually the center of $\mgoz$, generated by a Dehn twist about the boundary curve.
 Now $\autpisg$ is the quotient of $\mgoz$ by $\zz$.  In general,
 the quotient of a linear group is not necessarily linear, but the extra
 information we have about the kernel in this case will allow us to draw the desired conclusion.
  The following two theorems are proved in \cite{we}.  Note that the term
 ``closed'' refers to the Zariski topology.

 \begin{theorem}\label{chevlemma}
 Let $G$ be a linear group and $H$ a closed normal subgroup of $G$.  Then
 $G/H$ is also linear.  
 \end{theorem}

 \begin{theorem}\label{centersareclosed}
 The centralizer of any subset of a linear group is closed.
 \end{theorem}

 \proof[Proof of Lemma~\ref{alexlemma}]
 Since $\zz$ is the center of $\mgoz$, it is normal and also closed by
 the above.  Thus we can apply Theorem~\ref{chevlemma} to the surjection
 given in Sequence~\ref{surjection}, and Lemma~\ref{alexlemma} follows
 directly. \endpf

 We are now ready to prove the claim.

\proof[Proof of Claim~\ref{keyclaim}]
 Suppose that $\mgzo$ contains a poison subgroup.  Then by the isomorphism
 of Lemma~\ref{isolemma}, $\autpisg$ also contains a poison subgroup.
  Then $\autpisg$ is not linear by Lemma~\ref{hessamlemma}.  Now by Lemma~\ref{alexlemma},
 $\mgoz$ is also not linear. 
 The claim follows from the fact that $\mgoz$ is a subgroup of $\mgplusk$, for $k \geq 1$. Although this fact is well-known, for the sake of completeness we include a proof as follows.
Consider $S_{g,1,0}$ as a subsurface of $S_{g+k,0,0}$.
Let $h$ be the homomorphism from $\mgoz$ to $\mgplusk$ defined by extension to the identity on $S_{g+k,0,0} \setminus S_{g,1,0}$. Let $f \in \textrm{ker}(h)$ such that $f \neq id$. The mapping class $h(f)$ of $S_{g+k,0,0}$ keeps the subsurface $S_{g,1,0}$ invariant up to isotopy. According to Section 7.5 in \cite{I}, $h(f)$ induces a well defined mapping class in $\pi_0(\textrm{Diff}(S_{g,1,0}))$ (the group of homeomorphisms of $S_{g,1,0}$ up to isotopy not necessarily fixing $\partial S_{g,1,0}$). But since $h(f)=id$ and by the definition of $h$, this implies that $f$ induces the identity in $\pi_0(\textrm{Diff}(S_{g,1,0}))$, which implies that $f$ could only be a non-trivial power of a Dehn twist in the $\partial S_{g,1,0}$.  Then by definition, $h(f)$ will also be a non-trivial power of a Dehn twist, which is a contradiction. 
 \endpf

 \begin{remark}\label{myremark}
 {\rm We have defined poison subgroups in the context of $\mgzo$ and also in
 the context of automorphism groups, but the definition also makes sense
 in the context of any group action on another group.  Thus one could use this as a general approach to the linearity question for any such group.}  \end{remark}

 \section{Poison subgroups cannot be embedded in $\calm_{g,0,1}$}

Our strategy for proving this result will be to decompose the surface $S= S_{g,0,0}$ into subsurfaces in a particular way.  We then use the machinery of graphs of groups (described in detail in \cite{b}) to analyze the action of the generators of a poison subgroup of $\mgzo$ on the elements $x_1, x_2, x_3 \in \pi_1 (S)$.  
After completion of the proof of Theorem A, we discovered that similar methods involving graphs of groups and normal forms were used by Levitt and Vogtmann in \cite{lv} to give an algorithm for the Whitehead problem for surface groups.  There is a major difference, however, in that we are not given the curves $x_1, x_2,$ and $x_3$, and hence we cannot apply their algorithm directly, nor would our proof be significantly shortened by direct reference to their results.  Thus we have kept the proof of Theorem A in its original form for the sake of self-containment.  We have, however, found it useful to adopt their methods for the decomposition of the surface $S$. 

Throughout this section assume that $g \ge 2$, since Theorem A is clear when $g \le 1$. Fix a point $* \in S$, and identify $S_{g,0,1}$ with
 $(S,*)$. We
 use the point $*$ as the base point for the fundamental group of $S$.
 Let 
 $\langle \phi_1,\phi_2 \rangle$ be a poison subgroup in $\mgzo$. Then
 there are
 elements
  $x_1,x_2,x_3 \in
  \pi_1(S, *)$ such that $\langle x_1, x_2, x_3 \rangle  \cong F_3$ and

 \begin{enumerate}
   
 \item $\phi_i (x_j) = x_j$ \hskip .15in $ i, j = 1, 2$, and
  \item $\phi_i (x_3) = x_3 x_i$ \hskip .15in $i = 1, 2$.
  \end{enumerate}

In what follows, we will choose appropriate representatives for  
$\phi_i$ and $x_j$ (denoted by the same names by abuse of notation) such that, among other things, 
a power of $\phi_i$ fixes a regular neighborhood of $x_j$ pointwise. 
To this end our main tool will be the following result of Hass and Scott \cite{hs}. For $y_1,y_2 \in \pi_1(S,*)$, let
$$\textrm{Stab}(y_1,y_2)=\{\phi \in \mgzo \ | \ \phi(y_i)=y_i,  \ i=1,2 \}.$$

\begin{lemma}\label{hassy}
Let $y_1,y_2$ be distinct elements of $\pi_1(S,*)$, which are not proper powers. Then there exists a  representative of $y_i$ (denoted by $\tilde y_i$) and a subsurface $\cala$ formed by a regular neighborhood $N$ of $\tilde y_1 \cup \tilde y_2$ together with all disk components of $S \setminus N$, such that,
for any 
$ \phi \in {\rm Stab}(y_1,y_2)$, $\phi$ has a representative homeomorphism $\tilde \phi$ such that $\tilde \phi(\cala)=\cala$. 
\end{lemma}

This lemma follows from Theorem 2.1 in \cite{hs} together with the discussion in the beginning of page 32 in the same paper. For further details see Section 2.1 in \cite {lv}. 

\begin{remark}\label{uniqueona} {\rm Notice that, in Lemma~\ref{hassy}, if $ \phi \in \textrm{Stab}(y_1,y_2)$, the map $\phi$ induces a unique mapping class in $\pi_0(\textrm{Diff}(\cala,*))$ (see Section 7.5 in \cite{I}).}
\end{remark}

Since it is possible that $x_1$ and $x_2$ are proper powers, we need the following well-known lemma, adapted from \cite{lv}.

\begin{lemma}\label{properpowers}
Given a nontrivial element $x \in \pi_1(S,*)$ , 
there exists a unique $y \in \pi_1(S,*)$ and a unique $t\ge 1$ such that $y$ is not a proper power and $x = y^t$.
\end{lemma}

\proof  A proof is given in \cite{lv} (Lemma 2.3).  Though we will not give details, we note that it is also possible to prove this lemma by elementary hyperbolic geometry, using the discrete action of $\pi_1(S,*)$ on the upper half plane by hyperbolic isometries. \endpf

\begin{corollary}\label{uniqueroots}
Let $z_1,z_2 \in \pi_1(S,*)$ be such that $z_1^N=z_2^N$ for some $N \ge 1$. Then $z_1=z_2$.\end{corollary}

\proof Using Lemma~\ref{properpowers} let $y_i^{t_i}=z_i$ such that $y_i$ is not a proper power and $t_i \ge 1$, for $i=1,2$. Let  $x=y_1^{t_1N}=y_2^{t_2N}$. By the uniqueness guaranteed by Lemma~\ref{properpowers}, we have $y_1=y_2$ and $t_1N=t_2N$. Hence $z_1=z_2$, as desired.
\endpf 

Using Lemma~\ref{properpowers}, we can choose elements $y_i$ which are not proper powers and $t_i \ge 1$ such that $x_i = y_i^{t_i}$ for $i = 1,2$.  Then we know that $\phi_i (y_j^{t_j}) = y_j^{t_j}$, which implies that $\phi_i(y_j) =y_j$, by Corollary~\ref{uniqueroots}.
Notice that $y_1$ and $y_2$ are distinct since $\langle x_1, x_2 \rangle \cong F_2$. We choose $\tilde y_i$ and $\cala$ according to Lemma~\ref{hassy}.
Let $\pi_0(\textrm{Diff}(S,\cala))$ be the subgroup of $\mgzo$ consisting of mapping classes which have a representative keeping $\cala$ fixed pointwise.  We now adapt Lemma 3.1 of \cite{lv} to our purposes, and repeat their argument nearly verbatim.

\begin{lemma}\label{AC}
The subgroup  $\pi_0({\rm Diff}(S,\cala))$  has finite index in ${\rm Stab}(y_1,y_2)$. \end{lemma}

\proof
 First note that $\cala$ is not an annulus, since $x_1$ and $x_2$ generate a free group.
Using Lemma~\ref{hassy} (and noting Remark~\ref{uniqueona}), we can define a map $\rho$ from $ \textrm{Stab}(y_1, y_2)$ to $ \pi_0(\textrm{Diff}(\cala,*))$.
 Now we claim that the image of $\rho$ is finite.  To see this, let $k$ be 
any positive integer.  Let $T_k$ denote the set of homotopy classes of simple
 closed curves in $\cala$ whose intersection number with $y_1$ and
 $y_2$
 is at most $k$.  Then $T_k$ is finite, since $\cala \setminus (\tilde y_1
 \cup \tilde y_2)$ is composed entirely of disks and annuli.  Any map $\phi \in \textrm{Stab}(y_1, y_2)$ will preserve the intersection number of a curve with $y_1$ and $y_2$, and hence $\textrm{Stab}(y_1, y_2)$ acts on the set $T_k$.  Now choose a finite set $W$ of simple closed curves in $\cala$ whose image completely determines an element of $\pi_0(\textrm{Diff}(\cala,*))$. Let $k$ be bigger than the intersection number of any element in $W$ with $y_1$ and $ y_2$.  Thus the class of $\phi$ restricted to $\cala$ in $\pi_0(\textrm{Diff}(\cala, *))$ is completely determined by the action of $\phi$ on $T_k$.  But the set of
 permutations of $T_k$ is finite, and hence the image of $\textrm{Stab}(y_1,y_2)$ under $\rho$ is finite. 

Now let $\iota: 
\pi_0(\textrm{Diff}(\cala,*)) \to \textrm{Out}(\pi_1(\cala,*))$ be the natural homomorphism. The image of $\iota \circ \rho$ is also finite by the above argument.
Now any element $\phi \in \textrm{ker}(\iota \circ \rho)$ induces an inner automorphism on $\pi_1(\cala,*)$, i.e., $\phi(z)=czc^{-1}$. The element $c$ has to commute with both $y_1$ and $y_2$, which implies that $c$ has to be a power of both $y_1$ and $y_2$  since the centralizer of an element in a surface group is cyclic (this is an exercise in elementary hyperbolic geometry), and $y_1$ and $y_2$ are not proper powers. But this implies that $c=1$  since $x_1$ and $x_2$ generate a free group.
Hence $\phi$ induces the identity on $\pi_1(\cala,*)$. Picking a set of simple generators for $\pi_1(\cala,*)$, one can use an isotopy of the surface to make sure that $\phi$ keeps them fixed pointwise, by \cite{ep}.
Then one can further isotope $\phi$ to make sure $\phi$ keeps $\cala$ invariant pointwise by Alexander's lemma \cite{rolf}. Hence $\textrm{ker}(\iota \circ \rho)$ 
is contained in $\pi_0(\textrm{Diff}(S, \cala))$, which proves the lemma.
 \endpf

\begin{prop}\label{please}
There exists an integer $M$ such that $\phi_i^M$ fixes $\cala$ pointwise (up to isotopy).
\end{prop}

\proof  We know $\phi_i \in \textrm{Stab}(y_1,y_2)$ for $i = 1,2$. Hence by Lemma~\ref{AC}, there is an integer $M_{i} \ge 0$ such that $\phi_i^{M_i}\! \in\! \pi_0(\textrm{Diff}(S,\cala))$. Letting $M=\textrm{LCM}(M_1,M_2)$, we have $\phi_i^M \in\pi_0(\textrm{Diff}(S,\cala))$ for $i=1,2$. \endpf

From this point on, we assume that we are working with the particular representative of $\phi_i^M$ which fixes $\cala$ pointwise.

\begin{figure}\small
  $$
  \setlength{\unitlength}{0.05in}
  \begin{picture}(0,0)(0,11)
   \put(1,10){$\calb_3$}
   \put(17,25){$b_3$}
  \put(35,27){$\g_{3,1}$}
  \put(32,37){$\g_{3,2}$}
 \put(90,13){$\calb_2$}
 \put(76,33){$b_2$}
 \put(61,33){$\g_{2,1}$}
 \put(55,40){$e_{2,1}$}
 \put(47,40){$e_{3,1}$}
 \put(40,49.5){$e_{3,2}$}
 \put(51,53){$*$}
 \put(48.5,86){$\cala$}
 \put(65,71){$\g_{1,1}$}
 \put(53,61){$e_{1,1}$}
 \put(76,73){$b_1$}
 \put(90,63){$\calb_1$}
  \put(1,69){$\calb_r$}
   \put(17,80){$b_r$}
   \end{picture}
   \includegraphics[width=5in]{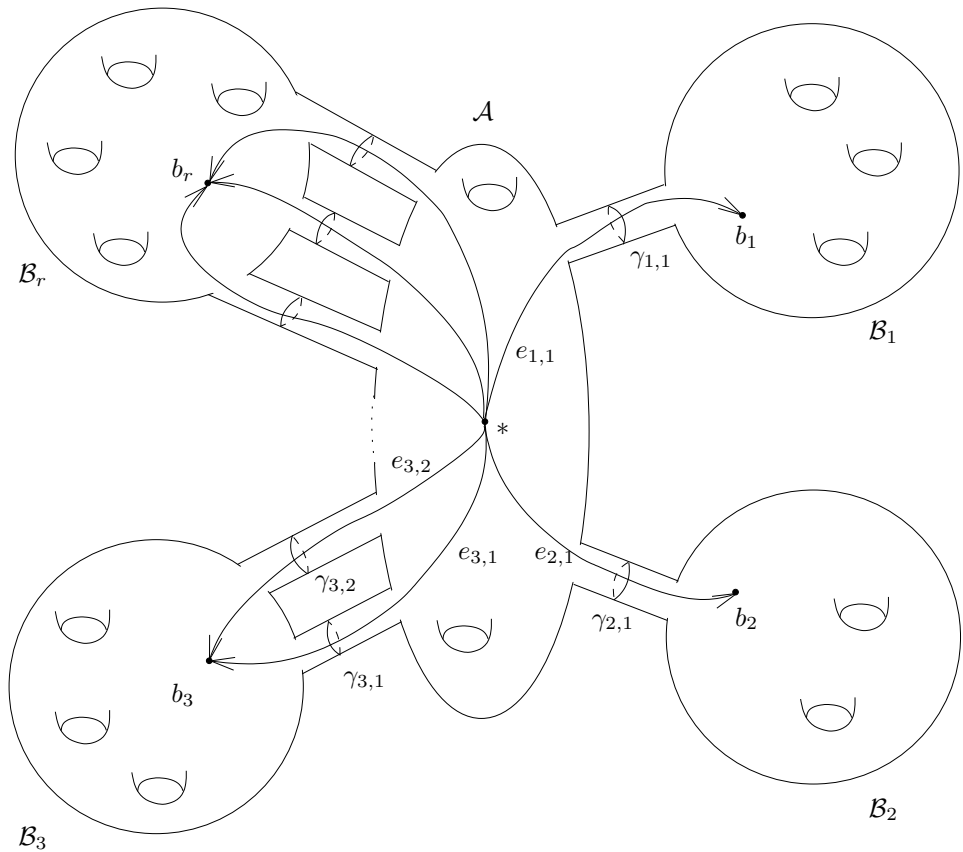}
  $$
 \caption{The decomposition of the surface $S$}
     \label{surface}
     \end{figure}

 \setlength{\unitlength}{0.05in}

 \begin{figure}\small
  
 $$\begin{picture}(0,0)(0,11)
   \put(-14,-7){$e_{j,k}$}
   \put(9,-8){$\g_{j,k}$}
  \put(-2,-10){$c_{j,k}$}
  \put(-2,-2){$e''_{j,k}$}
  \put(-2,-15){$e'_{j,k}$}
  \put(2,-22){$*$}
  \put(1,9){$b_j$}
   \end{picture}$$

  {\centerline{ \includegraphics[width=1in]{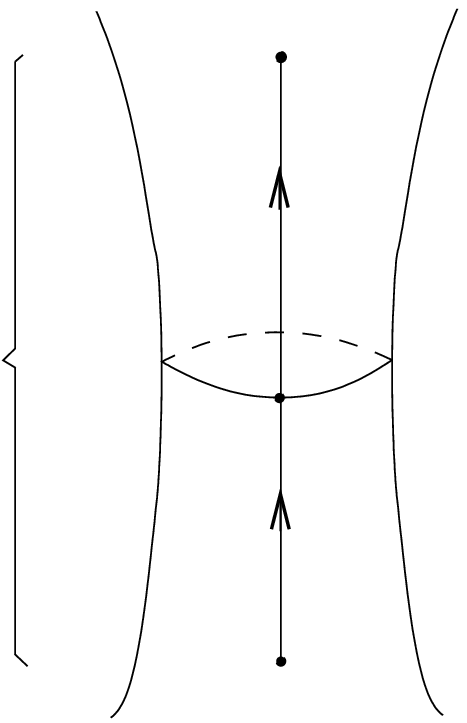}}}
     \caption{The subarcs of $e_{j,k}$}
     \label{fig2}
     \end{figure}

 Let $\calb_1,\cdots,\calb_r$
  be the respective 
  closures of  each
 component of
  $S \setminus \cala$.  Each component is $\calb_j$
  attached to $\cala$ along one or more circles.
 Hence $\cala \cap \calb_j$
  consists of $n_j \ge 1$
  circles, which we denote by $\g_{j,1},
   \cdots,\g_{j,{n_j}}$.

 In what follows we will use this decomposition of $S$ into the subsurfaces
 $\cala,\calb_j$ to construct a graph of groups $\calg$ whose
 fundamental group
 will give a decomposition of $\pi_1(S,*)$. To that end, we introduce
 some
 notation.

  For an oriented arc $e$ let $\st(e)$ and $\en(e)$ be the starting and
  ending points of the
  arc $e$, respectively. Also, let $\bar e$ be the same arc with the opposite
  orientation. In the
  following discussion, let the pair of indices $j,k$
 be such that $1 \le j \le r$, and $1\le k
 \le n_j$.

  Choose base points  $b_j \in \calb_j$. Notice
 that $\phi_i^M$ fixes each $\calb_j$ setwise. Hence we further isotope
 $\phi_i^M$ so that it fixes $b_j$, for $i=1,2$.
 See Figure~\ref{surface}.

 Choose oriented arcs $e_{j,k}$ connecting
 $*$ to $b_j$ for $1\le
  j \le r$ and $1 \le
  k \le n_j$. Choose each arc $e_{i,j}$  such that it intersects $\g_{j,k}$
 exactly once, and does not
  intersect any other $\g$'s. Moreover, we make the choices in such a
  way that if $(j,k) \ne (j',k')$, then
  $e_{j,k}$ and $e_{j',k'}$ do not intersect except possibly at the endpoints.
 Let $c_{j,k}$
 be the point of intersection of $e_{j,k}$ with
  $\g_{j,k}$.
   Also, let $e_{j,k}'$ be the subarc of $e_{j,k}$ connecting $*$ to $c_{j,k}$,
  and let $e_{j,k}''$ be
  the subarc from $c_{j,k}$ to $b_j$. See Figure~\ref{fig2}.

  Let $G$ be the graph embedded in $S$ with vertices $*,b_1,\cdots,b_r$
  and
  geometric edges  $e_{j,k}$ as above. As a technical point, the arcs with the opposite
  orientation 
  $\bar e_{j,k}$  are also considered edges of the graph $G$ but not drawn
  separately.

  We use the graph $G$ to construct a graph of groups. To each
  vertex of $G$ we 
  assign the fundamental group of the subsurface in which it is located,
  namely, to $*$ we
  assign $A=\pi_1(\cala,*)$, to
  $b_j$ we assign $B_j=\pi_1(\calb_j,b_j)$.
  To
  each edge
  $e_{j,k}$ we assign
  $\G_{e_{j,k}}=\pi_1(\g_{j,k}, c_{j,k}) \cong \zz$. Also, let $\G_{\bar e_{j,k}}=
  \G_{e_{j,k}}$. We also
  have  natural injections of the edge groups into the adjoining vertex
 groups
   as follows: for any $e_{j,k}$, since $\st(e_{j,k})=*$, the
  vertex group for 
  $\st(e_{j,k})$ is
  $A$. We have $\al_{e_{j,k}} :\G_{e_{j,k}} \to A$  defined by  $\al_{e_{j,k}}(x)=e'_{j,k}
  x \bar e'_{j,k}$.
  Corresponding to $\en(e_{j,k})$, we have $\bar \al_{e_{j,k}} :\G_{e_{j,k}}
 \to B_j$ which is defined by  $\bar
  \al_{e_{j,k}}(x)=\bar e''_{j,k}x
   e''_{j,k}$. For the edges $\bar e_{j,k}$ set $\al_{\bar e_{j,k}}=\bar \al_{e_{j,k}}$
 and  $\bar
 \al_{\bar
  e_{j,k}}= \al_{e_{j,k}}$.

   Let
  $\calg$ be the graph of groups constructed by the above data. By the
 generalized
  Van Kampen theorem, 
  $\pi_1(S,*)$ is isomorphic to the fundamental group of the graph of
 groups
  $\pi_1(\calg,*)$.

  To understand the elements of $\pi_1(\calg,*)$, we quote some definitions
  from \cite{b}. A
  {\it{loop}} based at $*$ in
  $\calg$ is a sequence
  $$t=(g_0,\e_1,g_1,\cdots,\e_n,g_n)$$ where $\e_i$ are edges of $G$ and
  $(\e_1,\cdots,\e_n)$
  is a loop in $G$ with $\st(\e_1)=*$ and $\en(\e_n)=*$. Also, $g_0$ and
 $g_n$ are in  
  $A$, and for
  $0 < i < n$, each
  $g_i$ is in the group assigned to $\en(\e_{i})=\st(\e_{i+1})$.
 A loop $t$  in $\calg$ is {\it{reduced}} if either $n=0$ and $g_0 \ne
  1$, or $n >0$ and whenever
  $\e_{i+1}=\bar \e_i$, we have $g_i
  \notin
  \al_{\bar \e_i}(\G_{\e_i})$. Geometrically, one can think of
  $t$ as a loop in $S$,
  with $g_i$ being loops in respective subsurfaces, and $\e_i$ as arcs
  connecting these loops. From this point of view, a reduced loop 
  on $S$ does not ``travel'' to a component $\calb_j$  unnecessarily.

   By \cite{b}, any non-trivial element of $\pi_1(\calg,*)$ can be written
  as $|t|=g_0\e_1g_1\cdots\e_ng_n$,
  where $t$ is a reduced loop as above. 

 \begin{remark}
{\rm
 The reduced loop representing 1 is the empty sequence.
}\end{remark}
\begin{remark} \label{nonreduced}
{\rm
  A non-reduced loop can be made into a reduced loop
 which represents the same element in $\pi_1(\calg, *)$
 by the process of {\it{combing}}. Namely, if a loop $t$ of length $n>1$
 is not reduced, it has a subsequence of the form 
 $(g_{i-1},\e_i, \al_{\bar \e_i}(h_i),\bar \e_i)$. One can replace this
 subsequence with $(g_{i-1} \al_{\e_i}(h_i))$. This process reduces the
 length, so after
 finitely many steps one arrives at a reduced loop. 
 }
   \end{remark}
   
 The following
  theorem is proved in \cite{b}.

  \begin{theorem}\label{uniqueness} Let $t=(g_0,\e_1,g_1,\cdots,\e_n,g_n)$
  and
  $t'=(g'_0,\e'_1,g'_1,\cdots,\e'_m,g'_m)$ be two reduced loops such that
  $|t|=|t'|$ in $\pi_1(\calg)$. Then $n=m$, $\e_i=\e'_i$ for $1 \le i
 \le
  n$,  and there exist
  $h_i \in \G_{\e_i}$ such that
\begin{enumerate}

\item  $g'_0=g_0  \ \al_{\e_i}(h_1)\inv$,

\item  $g'_i=\al_{\bar \e_{i}}(h_{i}) \ g_i  \ \al_{\e_{i+1}}(h_{i+1})\inv$,

\item  $g'_n=\al_{\bar \e_{n}}(h_{n}) \ g_n. $

\end{enumerate}

   \end{theorem}
   
 Notice that in the above theorem the elements of the form $\al_{\e}(h)$
 come from the circles $\g_{j,k}$.

\proof[Proof of Theorem A] Suppose $ \lan \phi_1,\phi_2 \ran \le \mgzo$  is a poison subgroup
 with
 respect
  to $x_1,x_2,x_3 \in
  \pi_1(S_{g,0,0},*)$. We construct the graph of groups $\calg$ as above,
  with $\pi_1(\calg, *) \cong
  \pi_1(S_{g,0,0},*)$. In the following we will identify these two groups.  

 By Proposition~\ref{please}, we can choose an integer $M$ such that $\phi_i^M$ fixes $\cala$ pointwise.  Since $\phi_i^M$ also sends each $\calb_j$
 to itself fixing the base points, we can see that
  $\phi_i^M(e_{j,k})=e_{j,k}p_{j,k}$ where
  $p_{j,k} \in B_j$. Similarly
  $\phi_i^M(\bar e_{j,k})= p_{j,k}\inv \bar e_{j,k}$.

 We will now simplify notation a bit by letting $\phi$ stand for $\phi_1^M$. 
 Let $x_3=|t|$ where $t$ is the reduced loop $t=(g_0,\e_1,g_1,\cdots,\e_{2n},g_{2n}).$
  Notice that since the graph $G$ is ``star-shaped'', the length of the
  loop must be even.
  Therefore
  $$\phi(x_3)=\hbox to 2.25in{\hss}$$
$$|(g_0, \e_1, p_1 \phi(g_1) p_2\inv\!, \e_2,g_2,\e_3,  p_3
  \phi(g_3) p_4\inv\!,
  \e_4,\cdots,p_{2n-1}\phi(g_{2n-1})p_{2n}\inv\!,\e_{2n},g_{2n})|$$ (each
  $p_i$ is in the group
  which makes this a well-defined path). Now by the condition $\phi_1(x_3)=x_3x_1$, which implies that $\phi(x_3) = x_3 x_1^M$,
 we get the equality
  $$|(g_0, \e_1, p_1 \phi(g_1) p_2\inv\!, \e_2,g_2,\e_3,  p_3 \phi(g_3)
 p_4\inv\!,
  \e_4,\cdots,p_{2n-1}\phi(g_{2n-1})p_{2n}\inv\!,\e_{2n},g_{2n})|$$
  $$= |(g_0,\e_1,g_1,\cdots,\e_{2n},g_{2n}x_1^M)|.$$
 Let $t'$ and $t''$ be the paths appearing on the left and right hand sides of the above equation respectively.
  Since the path $t$ is
  reduced, so is $t''$. If $t'$ is not reduced, by Remark~\ref{nonreduced}
 we can
 comb it to a reduced path $t'_{\rm{red}}$. By the equality
  and Theorem~\ref{uniqueness},
  $t'_{\rm{red}}$ must have the same length as $t''$, which means $t'$
 was reduced in
 the first place.
  Using Theorem~\ref{uniqueness} again, there is an $h_1 \in
  \G_{\e_{2n}}$ such that $g_{2n}x_1^M=\al_{\bar \e_{2n}}(h_1) \ g_{2n}$,
  i.e.,  $x_1^M=g_{2n}\inv
  \al_{\bar \e_{2n}}(h_1) \ g_{2n}$. Similarly, using $\phi_2$ in place
  of $\phi_1$, there
  exists an $h_2
  \in \G_{\e_{2n}}$ such that $x_2^M=g_{2n}\inv \al_{\bar \e_{2n}}(h_2)
 g_{2n}$. But $\G_{\e_{2n}} \cong \zz$, therefore $h_1,h_2$ commute, which
 implies $x_1^M,x_2^M$
  commute. This is a
  contradiction, since $\langle x_1,x_2 \rangle \cong F_2$. \endpf

\section{FP-groups do not embed in mapping class groups}
We begin by showing how to narrow our search for an FP-subgroup in a mapping class group.  

\begin{lemma} \label{reducetof2} Suppose that $\calm_{g,b,n}$ contains an
FP-subgroup. Then it contains an FP-subgroup $H$ which is isomorphic to a
quotient of $\calh(F_2)$. Moreover, the image of $F_2 \times \{1\}$ in $H$ is isomorphic to $F_2$. \end{lemma}

\proof Suppose $\calm_{g,b,n}$ contains an FP-subgroup. Hence there is a group $G$ and a homomorphism $\rho: \calh(G) \to \mgbn$ such that $\rho (G \times \{ 1 \})$ is not nilpotent-by-abelian-by-finite. Here $\rho( \calh(G)) \cong \calh(G)/\textrm{ker}(\rho)$ is the FP-subgroup of $\mgbn$. 
 By  Tits' alternative for mapping class groups (\cite{I} or
\cite{Mc}), $\rho(G \times \{1\})$  is either abelian-by-finite or contains a subgroup isomorphic to $F_2$. By assumption, the
latter holds. Let $x_1,x_2 \in G$ such that $\langle \rho(x_1,1), \rho(x_2,1) \rangle \cong F_2$. Then it is easily seen that for $G_1=\langle x_1,x_2 \rangle$, $\rho( \calh(G_1))$ is an FP-subgroup of $\mgbn$ and $G_1 \cong F_2$.  \endpf

We now recall the following definition from \cite{I}.
A mapping class $f$ is called {\it{pure}}
if there exists a set (possibly empty)
$\calc=\{c_1,\cdots,c_k\}$ of
 non-parallel, non-trivial, non-intersecting simple closed curves on the surface
such that:

\begin{enumerate}
\item The mapping class $f$ fixes each curve in $\calc$ up to isotopy.
\item The mapping class $f$ keeps each component of $S \setminus \calc$ invariant
up to isotopy.
\item  The restriction of $f$  to each component of $S \setminus \calc$  is either
the  identity or pseudo-Anosov. (Recall that the restriction of $f$ to a surface $U$ is pseudo-Anosov if and only if for any non-trivial simple closed curve $c$ in $U$ not isotopic to $\partial U$ and for any $N>0$, $f^N(c)$ is not isotopic to $c$.)

\end{enumerate}

For an integer $m$, let  $H_1(S, \zz/m\zz)$ be the first homology group of $S$ with coefficients in $\zz/m\zz$. We have an action of $\mgbn$ on  $H_1(S,\zz/m\zz)$, which defines a natural homomorphism $\mgbn \to {\rm{Aut}}( H_1(S,\zz/m\zz))$. The following theorem is due to Ivanov (\cite{I}, 1.8).

\begin{theorem} \label{niki}
For any integer $m \ge 3$, the group $$\Gamma_m= {\rm{ker}}(\mgbn \to {\rm{Aut}}( H_1(S,\zz/m\zz)))$$ is a normal subgroup of finite index in  $\mgbn$ consisting only of pure elements.
\end{theorem}

In the following discussion we will only need one such subgroup, so we set $m=3$ for simplicity. Any value $m\ge 3 $ would work as well. 

The reader should note that in the following theorem, the generators $\phi_i$, $\alpha_j$, and $t$ do not have precisely the same meaning as in Section 2.

\begin{theorem}\label{pure} Assume
$\calm_{g,b,n}$ contains an FP-subgroup. Then there exists an FP-subgroup of the
form
$H=\langle
\phi_1,\phi_2,\al_1,\al_2, t \rangle$ such that $\phi_1$, $\phi_2$, $\al_1$ and
$\al_2$ are in $\Gamma_3$ (in particular they are pure), and 
\begin{enumerate}
\item $\langle \phi_1 ,\phi_2 \rangle \cong F_2$, 
\item $\al_i$ commutes with $\phi_j$,
\item $t  (\phi_i \al_i) t\inv=\al_i$.
\end{enumerate}
\end{theorem}

\proof Let  $H$ be an FP-subgroup of the form $\rho(\calh(F_2))$ as in Lemma~\ref{reducetof2}, where $F_2=\langle x_1,x_2\rangle$.
Let $\al_i=\rho(1,x_i)$ and $\phi_i=\rho(x_i,1)$. By abuse of notation, we denote $\rho(t)$ by $t$.
Then  $H=\langle \phi_1,\phi_2,\al_1,\al_2, t \rangle$ is an FP-subgroup
satisfying (1) - (3) above, by definition of an FP-subgroup and Lemma~\ref{reducetof2}. Using Theorem~\ref{niki}, $\Gamma_3$ is a normal subgroup of $\mgbn$ of finite index. Let $N=[\mgbn:\Gamma_3]$.
Then $\al_i^N,\phi_i^N \in \Gamma_3$ are pure, and $\langle \phi_1^N,\phi_2^N \rangle \cong F_2$.
Replacing each of
$\al_i,\phi_j$ with their $N$th powers and keeping the same $t$, we get an FP-subgroup satisfying the conditions of the theorem.  \endpf

In the rest of this paper we assume that $\al_i,\phi_j$ and $t$ are maps as given in Theorem~\ref{pure}.
 
We can now exploit the machinery of pure mapping classes as developed in \cite{I}.
 For a pure mapping class $f$, one can always find a
 representative homeomorphism (which we will also denote by $f$) which fixes each curve in $\calc$ and
each component setwise.  Moreover, the mapping class $f$ induces well-defined mapping classes on
components of $S\setminus \calc$ (see Section 7.5 in \cite{I}). As an important technical point, for a component $T$ of $S \setminus \calc$, in order to get a 
well-defined mapping class $f|_T$ in the mapping class group of $T$, one should allow the isotopies in $T$ to move the points in
the components of $\partial T$ which are created as a result of cutting $S$ open.  Otherwise, an ambiguity results from combining $f|_T$ with a Dehn twist in a component of $\partial T$. In other words, when the surface is cut open along $\calc$, all the new boundary components which appear will be dealt with essentially as punctures.  The same remark holds when considering the mapping class group of a connected subsurface of $S$. 
 In what follows, the phrase ``up to isotopy'' will usually be dropped, but should be understood in any discussion of topological equivalence. 

In the above discussion, the collection $\calc$ corresponding to a pure mapping class $f$ may
not be canonical, but in fact one can always choose a canonical collection of isotopy classes of disjoint
simple closed curves, denoted by $\s(f)$, which we will define shortly.  For two 1-submanifolds $\calc_1$ and $\calc_2$ of $S$, let $$i(\calc_1,\calc_2)=\min \{| \calc'_1 \cap \calc'_2| \ | \ \calc'_i {\textrm{
is isotopic to }} \calc_i \}.$$  In other words, $i(\calc_1,\calc_2)$ is the geometric intersection number of $\calc_1$ and $\calc_2$. We then define $\s(f)$ by saying $c \in \s(f)$ if the two following conditions hold:
\begin{enumerate}
\item $f(c)=c.$
\item For any simple closed curve $\g$, if $i(\g,c) \ne 0$, then $f(\g) \ne \g$.
\end{enumerate} The collection $\s(f)$ is called the {\it essential reduction
system } for $f$. It is proved in \cite{I} (see Chapter 7) that $\s(f)$ is a finite
collection of disjoint simple closed curves, and $f$ restricted to each component
of $S\setminus \s(f)$ is either the identity or pseudo-Anosov.

 If $f \in \mgbn$ is not pure, then as discussed above there is some $N>0$ such
that
$f^N$ is pure.   Thus we can extend the definition of essential reduction systems by defining
$\s(f)$ to be equal to $\s(f^N)$. The notion of an essential reduction system  was originally defined
in \cite{BLM} for a mapping class, and was generalized in
\cite{I} to an arbitrary subgroup of $\mgbn$.  Note that $\s(f)$ is a topological
invariant of the mapping class
$f$.  We use this
notion to define an invariant for a pair of mapping classes in  $\mgbn$.

\begin{definition} {\rm For two mapping classes $f,h  \in \mgbn$, we let
$$i(f,h)=i(\s(f),\s(h)).$$ }
\end{definition}

Notice that this is invariant under simultaneous conjugacy:

\begin{prop}\label{invariance} For $t,f,h \in \mgbn$, $i(tft\inv,tht\inv)=i(f,h)$.
\end{prop}

\proof  First notice that $\s(tft\inv)=t(\s(f))$, for $f,t \in \mgbn$ (again see 
\cite{I}, Chapter 7).  Then we have that 
\begin{eqnarray*}
i(tft\inv,tht\inv) &=& i(\s(tft\inv),\s(tht\inv))\\
&=& i(t(\s(f)), t(\s(h)))\\
&=& i(\s(f),\s(h))\\
&=& i(f,h).
\end{eqnarray*}\vglue -10mm\endpf

The invariant $i(f,h)$ for $f,h \in \mgbn$ will be crucial in the proof of Theorem B.
We recall the following lemma, proved in \cite{I}.

\begin{lemma}[Ivanov]\label{x}
Let $f$ be a pure mapping class. If $X$ is a
subsurface or a simple closed curve on the surface such that $f^N(X)=X$ for some $N
\ge 1$, then 
$f(X)=X$. 
\end{lemma}

The following definition is also inspired by \cite{I}.

\begin{definition}{\rm  Let $f \in \mgbn$, and let
$T$ be the isotopy class of a connected subsurface of $S$. We say
$f$ keeps $T$  {\it{precisely invariant}} if $f(T)=T$ and if $f(c)\ne  c$ for each curve $c$ such that $i(c, \partial T) \neq 0$. }
\end{definition}

In particular we note that a pure mapping class $f \in \mgbn$ keeps all components of $S\setminus
\s(f)$ precisely invariant, by the basic property of $\s(f)$. Similarly, $f$ keeps each
regular neighborhood of $c \in
\s(f)$ precisely invariant.  We now develop a series of lemmas to prove Theorem B.

\begin{lemma} \label{comm} Let $f,\al$ be pure mapping classes in $\calm_{g,b,n}$
 such that
$\al f=f \al$. Let $T$ be a component of $S\setminus \s(f)$.
  Then we have

\begin{enumerate}
\item[{\rm (i)}]
 $\al(T)=T$, up to isotopy.

\item[{\rm (ii)}] $\al(c)=c$ for each $c \in \s(f)$.

\item[{\rm (iii)}]  $i(f,\al)=0$; i.e., $\s(f)$ and $\s(\al)$ can be isotoped off
each other.
\end{enumerate}

\end{lemma}

\proof  For any integer $N$, $\al^N$ commutes with $f$. This implies that 
$f(\al^N(T))$ $=\al^N(f(T))=\al^N(T).$ Suppose a simple closed curve $c$
intersects
$\partial
\al^N(T)$ non-trivially. Then $\al^{-N}(c)$ intersects
$\partial T$ non-trivially, and so 
$f(\al^{-N}(c)) \ne \al^{-N}(c)$, by assumption. Applying $\al^N$ to both sides, we
get
$f(c) \ne c$. Hence $f$ keeps $\al^N(T)$ precisely invariant. By the basic property of the essential reduction system, either $f|_T=id$ or $f|_T$ is pseudo-Anosov.

{\bf Case 1}\qua Assume $f|_T=id$. Since $  f|_{\al^N(T)}=(\al^N|_T) f|_{T} (\al^N|_T)^{-1},$  we have  $f|_{\al^N(T)}=id$ for all $N$. Notice that 
 $i(\partial \al^N(T),\partial T)=0$, since $f$ keeps $\al^N(T)$ precisely invariant for all $N$. Moreover, we claim that no component $c$  of $\partial \al^N(T)$ can be isotopic to a simple closed curve in $T$ which is not isotopic to a component of $\partial T$. Otherwise, one can find a simple closed curve $\g$ in $T$ such that $i(c,\g) \ne 0$. But $f(\g)=\g$, which contradicts the fact that $f$ keeps $\al^N(T)$ precisely invariant. Similarly one can show that no component of $\partial T$ can be isotopic to a simple closed curve in $\al^N(T)$ which is not isotopic to  $\partial \al^N(T)$.
 This shows that either $\al^N(T)=T$ or
$\al^N(T)$ can be isotoped off $T$. This in turn implies that the collection of
subsurfaces $\{ \al^N(T) \ |  \ N \in \zz \}$ is a collection of disjoint homeomorphic subsurfaces up to isotopy,  and hence it is a finite collection. This
shows that
$\al^N(T)=T$ for some $N$, and since $\al$ is pure, $\al(T)=T,$ up to isotopy, by Lemma~\ref{x}.

{\bf Case 2}\qua Let $f|_T$ be pseudo-Anosov. Again, since $  f|_{\al^N(T)}\!=\!(\al^N|_T\!) f|_{T} (\al^N|_T\!)^{-1},$  we have  $f|_{\al^N(T)}$ is pseudo-Anosov for all $N$. Also, notice that 
 $i(\partial \al^N(T),\partial T)=0$, since $f$ keeps $\al^N(T)$ precisely
 invariant for all $N$. Moreover, we claim that no component $c$  of
 $\partial \al^N(T)$ can be isotopic to a simple closed curve in $T$ 
which is not isotopic to a component of $\partial T$. Otherwise,
 since $c \in \partial \al^N(T)$ and $f$ is pure and pseudo-Anosov 
on $\al^N(T)$, we have $f(c)=c$. On the other hand, $c$ is in the interior of $T$ and $f$ is pseudo-Anosov on $T$, hence $f(c) \ne c$, which is a contradiction.  Similarly one can show that no component of 
$\partial T$ can be isotopic to a simple closed curve in $\al^N(T)$
 which is not isotopic to  $\partial \al^N(T)$.
 This shows that either $\al^N(T)=T$ or
$\al^N(T)$ can be isotoped off $T$. The rest of the argument is exactly as in 
Case 1. This proves (i).

To prove (ii), let $c \in \s(f)$. Let $T$ be component of $S\setminus \s(f)$ such that $c$ is a component of $\partial T$. Then $\al(T)=T$, by (i). This implies that $\al$ permutes the components of $\partial T$, which by Lemma~\ref{x} implies that $\al(c)=c$, proving (ii).

To prove (iii), let $c \in \s(f)$ and $\g \in \s(\al)$ such that $i(c,\g)>0$. Then by
definition of an essential reduction system, $\al(c)\ne c$, which contradicts
(ii). \endpf

Let $H=\langle \phi_1,\phi_2,\al_1,\al_2, t \rangle$ be an FP-subgroup of the type described in
Theorem~\ref{pure}. Notice that by Lemma~\ref{comm}(iii), $\s(\phi_i) \cup  \s(\al_j)$
is  collection of non-intersecting simple closed curves.
For
$i=1,2$, let
$C_i=\s(\al_i)
\cap
\s(\phi_i)$, 
$A_i=\s(\al_i) \setminus C_i$ and $D_i=\s(\phi_i) \setminus C_i$. Note that each
of $A_i,C_i$ or $D_i$ could be empty.

\begin{lemma}\label{atleast} For $i=1,2$, $A_i \cup D_i  \subset \s(\al_i \phi_i)$.
\end{lemma}

\proof  Without loss of generality, we prove $A_i \subset \s(\al_i \phi_i)$. Let $c
\in A_i$. Notice that by Lemma~\ref{comm}(ii), $\al_i(c)=\phi_i(c)=c$. If $c \notin
\s(\al_i\phi_i)$, by definition, there is a subsurface $U$ containing $c$ where
$U$ is a component of $S\setminus \s(\al_i \phi_i)$. Since $\al_i\phi_i|_U$ fixes $c$, it is not pseudo-Anosov and hence is the identity.
Similarly since
$c
\notin
\s(\phi_i)$, there is a subsurface $V$ containing
$c$ where $V$ is a component of $S\setminus \s(\phi_i)$ such that
$\phi_i|_V=id$. Therefore $\al_i|_{U \cap V}=id$. Since $c$ is not isotopic to
any component of 
$\partial U$ or
 $\partial V$, and $i(\partial U,\partial V)=0$, $c$ is not isotopic to any component of $\partial ( U \cap V)$. Then one can
find a simple closed curve $\g$ in $U \cap V$ such that $i(c,\g)>0$. But
$\al_i|_{U \cap V}=id$, so $\al_i(\g)=\g$, which contradicts the fact that $c \in \s(\al_i)$.
\endpf

\begin{lemma}\label{dunk}
$i(\phi_1,\phi_2)=0$.
\end{lemma}

\proof Recall that $\s(\al_i)=A_i \cup C_i$ and $\s(\phi_i)=C_i \cup D_i$. By
definition of essential reduction system and Lemma~\ref{comm}(ii), $i(\al_i, \phi_j)=0$
and so 
$$i(A_i,C_j)=i(A_i,D_j)=i(C_1,C_2)=i(C_i,D_j)=0,$$ for $i,j=1,2$. Therefore
$i(\al_1,\al_2)=i(A_1,A_2)$. Now by  Lemma~\ref{atleast},
\begin{eqnarray*}
i(\al_1\phi_1,\al_2\phi_2) &\ge&
i(A_1,A_2)+i(A_1,D_2)+i(D_1,A_2)+i(D_1,D_2)\\
&=&i(A_1,A_2)+i(D_1,D_2).
\end{eqnarray*}
 By part (3)
of Theorem~\ref{pure} and Proposition~\ref{invariance}, we have that 
\begin{eqnarray*}
i(A_1, A_2) &=& i(\al_1,\al_2) \\
&=& i(t (\phi_1 \al_1) t^{-1},t (\phi_2 \al_2) t^{-1}) \\
&=& i(\phi_1 \al_1, \phi_2 \al_2)\\
&\geq& i(A_1,A_2)+i(D_1,D_2).
\end{eqnarray*}
 Thus $i(D_1,D_2)=0$. Hence
\begin{eqnarray*} 
i(\phi_1,\phi_2)&=&i(\s(\phi_1), \s(\phi_2))\\
&=& i(C_1 \cup D_1 , C_2 \cup D_2)\\
&=& i(C_1,C_2)+i(C_1,D_2)+i(D_1,C_2)+i(D_1,D_2)\\
&=&0,
\end{eqnarray*} 
which proves the lemma. \endpf

For a connected subsurface $U$ of $S$, we define a subgroup $\Gamma_3(U)$ of the mapping class group of $U$ as follows:
$$\Gamma_3(U)=\{ f|_U \ | \ f \in \Gamma_3 \ {\rm{and}} \ f(U)=U\}.$$
Notice that all elements of $\Gamma_3(U)$ are pure. Also notice that if $\al_i$(respectively $\phi_i$) keeps $U$ invariant, then by Theorem~\ref{pure} we have $\al_i|_U \in \Gamma_3(U)$ (respectively $\phi_i|_U \in \Gamma_3(U)$).  The following lemma is proved in \cite{I} (Lemma 8.13).

\begin{lemma} \label{pacommpa}
Let $\Gamma$ be a subgroup of the mapping class group of a connected surface $U$ consisting of pure elements. If $f \in \Gamma$ is a pseudo-Anosov element, then its centralizer in $\Gamma$ is an infinite cyclic group generated by a pseudo-Anosov element.
\end{lemma}

\begin{corollary}\label{pacenter}
Let $\Gamma$ be a subgroup of the mapping class group of a connected surface $U$ consisting of pure elements. If $f,h \in \Gamma$ are  pseudo-Anosov elements, then either $f$ commutes with $h$ or their respective centralizers in $\Gamma$ intersect trivially.
\end{corollary}

\proof Let $C_{\Gamma}(f)$ denote the centralizer of $f$ in $\Gamma$. Suppose
there is an element $1 \ne \theta \in C_\Gamma(f) \cap C_\Gamma(h)$. Then $f,h \in C_\Gamma(\theta)$, which is cyclic by Lemma~\ref{pacommpa}, so $f$ commutes with $h$. \endpf

We are going to encounter the following  particular situation in different contexts, so we declare it a lemma:

\begin{lemma}\label{centralizer}
Let $U$ be a component of $S \setminus \s(\phi_i)$ for $i=1$ or $i=2$ such that $\Gamma_3(U)$ is non-trivial. Assume that $\al_i|_U=id$ and  $\phi_i(U)= U$ for $i=1,2$.
 Then the respective centralizers of $\phi_1|_U$ and $\phi_2|_U$ in $\Gamma_3(U)$ intersect non-trivially.
\end{lemma}

\proof 
Without loss of generality, let $U$ be a component of $S\setminus \s(\phi_1)$. Assume on the contrary that  the centralizers of $\phi_1|_U$ and $\phi_2|_U$
in the mapping class group of $U$ have only the identity map in common. This 
in particular implies that $\phi_i|_U \ne id$ for $i=1,2$. 
The map $\phi_1|_U$ is pseudo-Anosov, since $U$ is a component of $S \setminus \s(\phi_1)$. Consider the subsurface $t(U)$. By part (3)
of Theorem~\ref{pure}, we have
\begin{equation}\label{conj}
\al_i|_{t(U)}= (t|_{U})( \phi_i|_{U} \al_i|_{U})( t|_{U})^{-1} =
 ( t|_{U})( \phi_i|_{U})( t|_{U})^{-1}.
\end{equation}
This implies that $\al_i|_{t(U)}\ne id$ keeps $t(U)$ invariant, since it is conjugate to $\phi_i|_{U}$, for $i=1,2$. Moreover, $\al_1|_{t(U)}$ is pseudo-Anosov. This in particular implies that $t(U)$  is a component of $S \setminus \s(\al_1)$, and $t(U)$ can be
isotoped off $U$, since $\al_1|_U=id$. Moreover, by assumption and by (\ref{conj}), the centralizers of $\al_1|_{t(U)}$ and $\al_2|_{t(U)}$ intersect trivially in  $\Gamma_3(t(U))$.  By Lemma~\ref{comm}(i), $\phi_i$ keeps $t(U)$ invariant for $i=1,2$, since $\phi_i$ commutes with $\al_1$. Again, since $\phi_i|_{t(U)}$ commutes with $\al_j|_{t(U)}$ and by the assumption about the centralizers, we have $\phi_i|_{t(U)}=id$, for $i,j=1,2$.
Now we can prove the following statements for $N \ge 1$ simultaneously by induction on $N$:

\begin{enumerate}
\item $\al_i|_{t^N(U)}\ne id$ keeps $t^N(U)$ invariant, for $i=1,2$.  \item $\al_1|_{t^N(U)}$ is pseudo-Anosov (hence, $\phi_i$ keeps $t^N(U)$ invariant for $i=1,2$).
\item The respective centralizers of $\al_i|_{t^N(U)}$ in $\Gamma_3(t^N(U))$ intersect trivially, for $i=1,2$.
\item  $\phi_i|_{t^N(U)}=id$, for $i=1,2$. 
\end{enumerate}
We have already established all four statements for $N=1$. The passage from $N$ to $N+1$ follows similarly from the relation:
$$\al_i|_{t^{N+1}(U)}\!=\! (t|_{t^N(U)})( \phi_i|_{t^N(U)} \al_i|_{t^N(U)})( t|_{t^N(U)})^{-1}\!\! =\!
 ( t|_{t^N(U)})( \al_i|_{t^N(U)})( t|_{t^N(U)})^{-1}$$

The second statement above shows that 
$t^N(U)$ can be isotoped off $U$, since $\al_1|_U=id$. Therefore, $t^M(U)$  can be isotoped off $t^N(U)$
for all $M \ne N$. This is clearly a contradiction, since the Euler characteristic
of $S$ is finite. \endpf

\begin{lemma}\label{clashingpas} For $i = 1,2$,
 let $U_i$ be a component of $S\setminus \s(\phi_i)$ such that $\phi_i|_{U_i}$is pseudo-Anosov. Then either $U_1$ and $U_2$ are disjoint up to isotopy, or $U_1$ is isotopic to $U_2$. 
\end{lemma}

\proof  First we show that if $U_1$ and $U_2$ are not disjoint, then either $U_1 \subseteq U_2$ or $U_2 \subseteq U_1$.  Suppose $U_1 \nsubseteq U_2$ and $U_2 \nsubseteq U_1$ but $U_1$ cannot be
isotoped off $U_2$. Throughout the proof, let $j,k \in \{1,2\}$ be arbitrary such that $j \neq k$. Since $i(\partial U_1,
\partial U_2)=0$, there is some component $c_j$ of $\partial U_j$ such that $c_j
\subset  U_{k}$ and $c_j$ is not isotopic to any component of $\partial U_{k}$.  
 By Lemma~\ref{comm}(i),
$\al_i$ keeps $U_1$ and $U_2$ invariant for $i = 1,2$. Since $\al_i \in \Gamma_3$, we have  $\al_i|_{U_j} \in \Gamma_3(U_j)$. Since $c_j$ is in the interior of $U_k$ and
$\al_i(c_j)=c_j$ by Lemma~\ref{comm}(ii),  this implies that $\al_i|_{U_k}$ is not pseudo-Anosov, hence by Lemma~\ref{pacommpa}, $\al_i|_{U_k}=id$  for $i,k=1,2$. 

Let
$U=U_1 \cup U_2$. 
At this point we apply a similar argument as in the proof of Lemma~\ref{centralizer}, as follows.  By the relation
\begin{equation}\label{conj1}
\al_i|_{t(U_i)}= (t|_{U_i})( \phi_i|_{U_i} \al_i|_{U_i})( t|_{U_i})^{-1} =
 ( t|_{U_i})( \phi_i|_{U_i})( t|_{U_i})^{-1},
\end{equation}
we see that $\al_i|_{t(U_i)}$  is pseudo-Anosov. This in particular implies that $t(U)=t(U_1)\cup t(U_2)$ can be isotoped off $U$, since $\al_i|_U=id$. Note that $t(U_i)$ is a component of $S \setminus \s(\al_i)$, so $\phi_j$ keeps $t(U_i)$ invariant for $i,j=1,2$, by Lemma~\ref{comm}(i). Since $\phi_i$ is pure, and $t(c_j)$ is a boundary component of $t(U_j)$, we have  $\phi_i(t(c_j))=t(c_j)$. By the choice of $c_j$ we know that $t(c_j)$ is in the interior of $t(U_k)$. By Lemma~\ref{pacommpa} and the fact that $\phi_i|_{t(U_k)} \in \Gamma_3(U_k)$, we have $\phi_i|_{t(U_k)}=id$, for $i,k=1,2$.
Now by induction on $N$ we can simultaneously prove the following statements for $N \ge 1$:

\begin{enumerate}
\item The map $\al_i|_{t^N(U_i)}$  is pseudo-Anosov, for $i=1,2$.
\item We have $\phi_i|_{t^N(U_j)}=id$, for $i,j=1,2$.
\end{enumerate}
We have already established these two statements for $N=1$.
The passage from $N$ to $N+1$ can be achieved by considering the conjugacy relation
\begin{equation}\label{conjN}
\al_i|_{t^{N+1}(U_i)}= t|_{t^N(U_i)} \ \phi_i|_{t^N(U_i)}  \ \al_i|_{t^N(U_i)} \ t|_{t^N(U_i)}^{-1} =
 t|_{t^N(U_i)} \  \al_i|_{t^N(U_i)} \  t|_{t^N(U_i)}^{-1}.
\end{equation}
This proves statement (1) above. Now  use Lemma~\ref{comm}(i) to see that $\phi_i$ keeps $t^{N+1}(U_j)$ invariant. This implies that $\phi_i|_{t^{N+1}(U_j)} \in \Gamma_3(t^{N+1}(U_j))$, and by Lemma~\ref{pacommpa}, we have statement (2).

In particular, statement (1) shows that
$t^N(U)$ can be isotoped off $U$ for all $N>1$, which is a contradiction as in Lemma~\ref{centralizer}. This proves that either $U_1 \subseteq U_2$ or $U_2 \subseteq U_1$, or $U_1$ and $U_2$ can be isotoped off each other.

Now without loss of generality, suppose that $U_1 \subseteq U_2$, but $U_1$ is not isotopic to $U_2$. Then there exists a component $c_1$ of $\partial U_1$ such that $c_1$ is not isotopic to a component of
$\partial U_{2}$. 
 By Lemma~\ref{comm}(i),
$\al_i$ keeps $U_1$ and $U_2$ invariant
 for $i = 1,2$. Also, by Lemma~\ref{comm}(ii), $\al_i(c_1)=c_1$,
 which implies $\al_i|_{U_2}=id$, by Lemma~\ref{pacommpa}. 
Again, using (\ref{conj1}) we get statement (1) for $N=1$.
 Hence $t(U_i)$ is a component of $S \setminus \s(\al_i)$. So $\phi_i$ keeps
 $U_j$ invariant. Thus $\phi_i(t(c_1))=t(c_1)$, which gives $\phi_i(U_2)=id$, by Lemma~\ref{pacommpa}. This proves statement (2) for $N=1$. The passage
 from $N$ to $N+1$ follows by using equation~(\ref{conjN}) above. Then again we have that $U_2$
 can be isotoped off $t^N(U_2)$ for all $N>1$, which is a contradiction. This proves 
that $U_1$ is isotopic to $U_2$. \endpf

\begin{lemma}\label{commutingpas} Let $U$ be a component of both $S\setminus
\s(\phi_1)$ and
$S\setminus \s(\phi_2)$ such that $\phi_i|_{U}$ is pseudo-Anosov for $i=1,2$. Then $\phi_1|_U$ commutes with $\phi_2|_U$.
\end{lemma}

\proof  If $\phi_1|_U$ and $\phi_2|_U$ do not commute, then their centralizers in $\Gamma_3(U)$ have trivial intersection by Corollary~\ref{pacenter}.  This implies that $\al_i|_U=id$, which contradicts Lemma~\ref{centralizer}.
 \endpf

We are finally ready to prove Theorem B.  

\proof[Proof of Theorem B] 
 Let $U$ be a component of $S\setminus \s(\phi_1)$ such that $\phi_1|_U$ is pseudo-Anosov. We first prove that $\phi_2|_U$ is either pseudo-Anosov or the identity. Suppose $\phi_2|_U$ is neither pseudo-Anosov nor the identity (in particular, $U$ is not a component of $S\setminus \s(\phi_2)$). Let $V_1,V_2, \cdots, V_s$ be components of 
$S\setminus \s(\phi_2)$, which cover $U$ up to isotopy. We can assume that the cover is minimal in the sense that none of the $V_k$ can be isotoped off $U$.
By Lemma~\ref{clashingpas}, $\phi_2|_{V_k}=id$ for all $1 \le k \le s$.  (This does not mean that $\phi_2|_U=id$, since $\phi_2$ may involve Dehn twists about boundary components of $V_k$.)  By Lemma~\ref{dunk},  $i(\partial U, \partial V_k)=0$ for all $1 \le k \le s$, which shows that $\phi_2$ keeps $U$ invariant.  Moreover, $\phi_2|_U$ is a non-trivial composition of Dehn twists about disjoint simple closed curves. Using Lemma~\ref{comm}(i), $\al_i$ keeps $U$ invariant. Since $\al_i|_U, \phi_j|_U  \in \Gamma_3(U)$ and $\al_i|_U$ commutes with $\phi_2|_U$, using Lemma~\ref{pacommpa} we see that $\al_i|_U$ cannot be pseudo-Anosov.
Moreover, $\al_i|_U$ commutes with $\phi_1|_U$ so $\al_i|_U=id$.
Now by Lemma~\ref{centralizer}, we get that the centralizers of $\phi_1|_U$ and $\phi_2|_U$ must intersect non-trivially. Lemma~\ref{pacommpa} then implies that $\phi_2|_U$ is either pseudo-Anosov or the identity, which is a contradiction.

We have proved that for a component $U$ of $S \setminus \s(\phi_1)$ where $\phi_1|_U$ is pseudo-Anosov,  $\phi_2|_U$ is either pseudo-Anosov or the identity. In the case that $\phi_2|_U$ is pseudo-Anosov, $\phi_1|_U$ and $\phi_2|_U$ commute by Lemma~\ref{commutingpas}. Similarly, for a component $V$ of $S \setminus \s(\phi_2)$ where $\phi_2|_V$ is pseudo-Anosov,  $\phi_1|_V$ is either a commuting pseudo-Anosov or the identity.

Let $S_1$ be the subsurface of $S$ which is the union of subsurfaces $T$ such that either $\phi_1|_T$ or $\phi_2|_T$ is pseudo-Anosov. We have proved that $\phi_1$ and $\phi_2$ both keep $S_1$ invariant, and $\phi_1|_{S_1}$ commutes with $\phi_2|_{S_1}$.

On $S_2=S \setminus S_1$ both $\phi_1$ and $\phi_2$ are compositions of Dehn twists about disjoint curves, by Lemma~\ref{dunk}. Hence $\phi_1|_{S_2}$ and $\phi_2|_{S_2}$ commute. 
We conclude that $\phi_1$ and $\phi_2$ commute,
contradicting part (1) of Theorem~\ref{pure}. This shows that FP-groups do not embed in $\mgbn$, as desired. \endpf

\Addresses

 \end{document}